\documentclass[10pt]{article}

\usepackage[english]{babel}
\usepackage{latexsym, graphicx, natbib}

\DeclareSymbolFont{AMSb}{U}{msb}{m}{n}
\DeclareSymbolFontAlphabet{\mathbb}{AMSb}

\begin{document}
\setcounter{page}{1}
\setcounter{section}{0}
\newcommand{\nc}{\newcommand}
\nc{\nt}{\newtheorem} \nt{defn}{Definition}[section] \nt{lem}{Lemma}[section]
\nt{pr}{Proposition}[section]
%\nt{th}{Theorem}
\nt{theorem}{Theorem}[section] \nt{cor}{Corollary}[section] \nt{ex}{Example}[section]
\nt{note}{Remark}[section] \nc{\bd}{\begin{defn}}
\nc{\ed}{\end{defn}} \nc{\blem}{\begin{lem}} \nc{\elem}{\end{lem}}
\nc{\bpr}{\begin{pr}} \nc{\epr}{\end{pr}} \nc{\bth}{\begin{theorem}}
\nc{\eth}{\end{theorem}} \nc{\bcor}{\begin{cor}}
\nc{\ecor}{\end{cor}} \nc{\bex}{\begin{ex}}  \nc{\eex}{\end{ex}}
\nc{\bnote}{\begin{note}}  \nc{\enote}{\end{note}} \nc{\prf}{{\bf
Proof} }
\renewcommand{\contentsname}{Table of Contents}
\nc{\eop}{\hfill $\Box$ \\ \\}
\nc{\disc}{\mathrm{Disc}}
\nc{\homo}{\mathrm{Hom}}
\nc{\vol}{\mathrm{vol}}
%\nc{\deg}{\mathrm{deg}}
\nc{\trace}{\mathrm{Tr}}
\nc{\hm}{\mathrm{Hom}}
\nc{\class}{\mathrm{Cl }}
\title{Riemann-Roch and Riemann-Hurwitz theorems for global fields}
\nc{\chara}{\mathrm{char}}
\author{Stella Anevski}
\date{}
\maketitle
\noindent
In this paper, we use counting theorems from the geometry of numbers to extend the Riemann-Roch theorem and the Riemann-Hurwitz formula to global fields of arbitrary characteristic. \\
 \\
In the first part of the paper (cf. \S \ref{gfatv}), we review some basic definitions and results from the theory of valuations on global fields. \\
 \\%We choose canonical representative from each equivalence class, and state the product formula. Define the integers ...$S_i$ determined by the Archimedean valuations. We also define the integers $e_Q$ and real numbers $r_Q$, determined by a finite separable extension of global fields.\\
In the following section (cf. \S \ref{dogfatrr}), we define a {\em divisor on a global field} ${\mathbb K}$ to be a formal linear combination of equivalence classes of valuations on ${\mathbb K}$ (cf. Definition \ref{divdef}). If the characteristic of ${\mathbb K}$ is zero, this definition coincides with the notion of an Arakelov divisor. If the characteristic of ${\mathbb K}$ is positive, the definition coincides with the ordinary notion of a divisor on the complete non-singular curve with function field ${\mathbb K}$. We then deviate from standard fare, by defining the degree of a divisor $D$ multiplicatively (cf. Definition \ref{divdef}). (The usual notion of the degree of $D$ is simply recovered by taking the logarithm to a suitable base.) \\
 \\For a divisor $D$, we define its set of {\em multiples} $H^{0}(D)$ in a standard manner (cf. Defintion \ref{muldef}). However, we cannot define $h^{0}(D)$ in the usual way, as the dimension of $H^{0}(D)$, since this would not make sense when the characteristic of ${\mathbb K}$ is zero. Instead, we declare that $h^{0}(D)$ is the cardinality of $H^{0}(D)$.\\
\\
With these definitions, we are in a position to formulate a Riemann-Roch type result for global fields (cf. Theorem \ref{rr}), using a counting theorem for metrized modules, proved by H. Gillet and C. Soul\'e (cf. \cite{GiSo}) in the characteristic zero case. The statement is that there exists a divisor $\omega_{\mathbb K}$ on the global field ${\mathbb K}$ such that
\begin{eqnarray*}
\frac{1}{C({\mathbb K})}\leq \frac{h^{0}(D)}{h^{0}(\omega_{\mathbb K}-D)}\cdot \frac{\sqrt{\deg \omega_{\mathbb K}}}{\deg D}\leq C({\mathbb K}),
\end{eqnarray*}
for any divisor $D$ on ${\mathbb K}$. \newpage
\noindent
In the inequalities above, $C({\mathbb K})$ denotes a quantity depending only on the archimedean valuations on ${\mathbb K}$. In particular, $C({\mathbb K})=1$ iff ${\mathbb K}$ has no archimedean valuations, in which case we recover the classical Riemann-Roch theorem by taking logarithms. We also use a theorem from \cite{Neu} to conclude that there exists a function $i$ on the set of divisors on ${\mathbb K}$ such that $i(\cdot)\to1$ as $\deg \cdot\to\infty$, and such that
\begin{eqnarray*}
\frac{h^{0}(D)}{i(D)}\cdot\frac{\sqrt{\deg \omega_{\mathbb K}}}{\deg D}=B({\mathbb K}),
\end{eqnarray*}
for any divisor $D$ on ${\mathbb K}$. Here $B({\mathbb K})$ denotes a quantity depending only on the archimedean valuations on ${\mathbb K}$, such that $B({\mathbb K})=1$ iff ${\mathbb K}$ has no archimedean valuations. Once again, we recover the classical Riemann-Roch theorem in this case.\\
 \\
In the final section (cf. \S \ref{atrht}), we associate to a finite separable extension of global fields ${\mathbb L}/{\mathbb K}$, a divisor $R_{{\mathbb L}/{\mathbb K}}$ on ${\mathbb L}$. For a special choice of ${\mathbb K}$ (cf. Definition \ref{knoll}), this gives us a "canonical divisor" $\omega_{\mathbb L}$ on ${\mathbb L}$ that can be used in Theorem \ref{rr}. Finally, we note that the divisors so defined satisfy a formula
\begin{eqnarray*}
\deg{\omega_{\mathbb L}}=\deg{\omega_{\mathbb K}}^{[{\mathbb L}:{\mathbb K}]}\cdot\deg{R_{{\mathbb L}/{\mathbb K}}},%\cdot \sqrt{2^{S_{\infty}(2)([{\mathbb L}:{\mathbb K}]-1)}}.
\end{eqnarray*}
which coincides with the classical Riemann-Hurwitz formula when the global fields have positive characteristic.\\
 \\
\section{Global fields and their valuations} \label{gfatv}
In this section, we introduce some notation, and give a brief review of the theory of valuations on global fields.
\bd
Let ${\mathbb K}$ be a field. A {\em valuation} on ${\mathbb K}$ is a mapping $\varphi:{\mathbb K}\to {\mathbb R}^{+}$, such that
\begin{eqnarray*}
(i) &&\mbox{$\varphi(\alpha)=0$ if and only if $\alpha=0$,}\\
(ii) &&\mbox{$\varphi(\alpha\beta)=\varphi(\alpha)\varphi(\beta)$,}\\
(iii) &&\mbox{$\varphi(\alpha+\beta)\leq \varphi(\alpha)+\varphi(\beta)$,}
\end{eqnarray*}
for all $\alpha,\beta\in {\mathbb K}$. A valuation $\varphi$ on ${\mathbb K}$ is said to be {\em archimedean} if
\begin{eqnarray*}
\varphi(\alpha+\beta)>\max(\varphi(\alpha),\varphi(\beta)),
\end{eqnarray*}
for some $\alpha,\beta\in {\mathbb K}$.
\eop
\ed
We define an equivalence relation $\sim$ on the set of valuations on ${\mathbb K}$, by declaring that $\varphi\sim \phi$ if and only if
\begin{eqnarray*}
\varphi(\alpha)<1 \Leftrightarrow \phi(\alpha)<1, \mbox{ for all $\alpha\in {\mathbb K}$.}
\end{eqnarray*}
%(This condition is clearly equivalent with the existence of a $c\in {\mathbb R}^{+}$ such that
%\begin{eqnarray*}
%\varphi(\alpha)=\phi(\alpha)^{c},
%\end{eqnarray*}
%for all $\alpha\in {\mathbb K}$.)
We denote by $\sum_{\mathbb K}$ the set of equivalence classes under this relation.\\
 \\
It is easily verified that whenever a valuation $\varphi$ is archimedean, $\varphi\sim\phi$ implies that the valuation $\phi$ is also archimedean. We say that an element $P\in\sum_{\mathbb K}$ is {\em archimedean} if $P$ contains an archimedean valuation, and we denote by $\sum_{\mathbb K}^{\infty}$ the set of archimedean elements in $\sum_{\mathbb K}$.  \\
 \\
\bd \label{knoll}
By a {\em global field} we mean either
\begin{eqnarray*}
- &&\mbox{a finite extension of the field ${\mathbb Q}$, or}\\
- &&\mbox{a finite extension of a field ${\mathbb F}_q(t)$ of rational functions in an}\\
&&\mbox{indeterminate $t$ over a finite field ${\mathbb F}_q$.}
\end{eqnarray*}
For a global field ${\mathbb K}$, we set
\begin{eqnarray*}
{\mathbb K}_0&=&\left\{ \begin{array}{lcl}
                                {\mathbb Q}, &&\mbox{ if $\chara({\mathbb K})=0$, }\\
                                {\mathbb F}_{q}(t), &&\mbox{ if $\chara({\mathbb K})>0$}.
                                         \end{array} \right .
\end{eqnarray*}
\eop
\ed
\bnote
If ${\mathbb K}$ is a global field of characteristic zero, the extension ${\mathbb K}/{\mathbb K}_0$ is separable (cf. {\em \cite{Lang}}, Corollary $6.12$ in \S V:$6$). For a global field ${\mathbb K}$ of positive characteristic, this is not always the case for an arbitrary choice of $t$. However, there is at least one choice of $t$ which makes the extension ${\mathbb K}/{\mathbb K}_0$ separable (cf. {\em \cite{Lang}}, Proposition $4.9$ in \S VIII:$4$), and in the sequel, we shall assume that such a choice is made in Definition \ref{knoll}.
\eop
\enote
%\bd
%A {\em place} of a global field ${\mathbb K}$ is an equivalence class of valuations on ${\mathbb K}$ under the relation $\sim$. A place is said to be {\em archimedean} if it contains an archimedean valuation.
%\eop
%\ed
%We denote the set of places of a global field ${\mathbb K}$ by $\sum_{\mathbb K}$, and the set of archimedean places of ${\mathbb K}$ by $\sum_{\mathbb K}^{\infty}$.
Let ${\mathbb K}$ be a global field, and choose a representative $\phi_P\in P$, for each $P\in \sum_{\mathbb K}$. Set
\begin{eqnarray*}
&&A_P=\{\alpha\in {\mathbb K}; \phi_P(\alpha)\leq 1\},\\
&&M_P=\{\alpha\in {\mathbb K}; \phi_P(\alpha)<1\}.
\end{eqnarray*}
Denote by $\widehat{{\mathbb K}}_P$ the completion of ${\mathbb K}$ with respect to $\phi_P$. We define a function $N:\sum_{\mathbb K}\to [1,\infty)\subset {\mathbb R}$, by letting
\begin{eqnarray*}
N(P)&=&\left\{ \begin{array}{lcl}
                                e^{\dim_{\mathbb R}\widehat{{\mathbb K}}_P}, &&\mbox{ if $P\in\sum_{\mathbb K}^{\infty}$, }\\
                                \#\left(A_P/M_P\right), &&\mbox{ if $P\in \sum_{\mathbb K}\setminus\sum_{\mathbb K}^{\infty}$}.
                                         \end{array} \right .
\end{eqnarray*}
Indeed, it is easily seen that both $\widehat{{\mathbb K}}_P$ and $A_P/M_P$ are independent of the choice of $\phi_{P}\in P$. Since ${\mathbb K}$ is global, the residue field $A_P/M_P$ is finite for all $P\in\sum_{\mathbb K}\setminus \sum_{\mathbb K}^{\infty}$ (cf. \cite{Rein}). For $P\in \sum_{\mathbb K}^{\infty}$, the completion $\widehat{{\mathbb K}}_P$ is either ${\mathbb R}$ or ${\mathbb C}$. Hence $N(P)<\infty$ for all $P\in\sum_{\mathbb K}$. \\
 \\
We define the integers $S_{i}({\mathbb K})=i\cdot \#\left(\{P\in\sum_{{\mathbb K}}^{\infty}; \log N(P)=i\}\right)$, $i\in\{1,2\}$, determined by the archimedean valuations on ${\mathbb K}$. 
\bnote
With this notation, $S_{1}({\mathbb K})=S_{2}({\mathbb K})=0$ if $\chara({\mathbb K})>0$. If $\chara({\mathbb K})=0$, then $S_{1}({\mathbb K})$ is the number of real embeddings of ${\mathbb K}$, and $S_{2}({\mathbb K})$ is the number of complex embeddings of ${\mathbb K}$.
\eop
\enote 
For $P\in\sum_{\mathbb K}\setminus\sum_{\mathbb K}^{\infty}$, we define the normalized valuation $\varphi_P\in P$ by requiring that
\begin{eqnarray*}
-\log_{N(P)}\varphi_{P}({\mathbb K})={\mathbb Z}\cup \{{\infty}\}.
\end{eqnarray*}
%Indeed, a valuation $\varphi$ belongs to $P$ if and only if there exists a $c\in {\mathbb R}^{+}$ such that
%\begin{eqnarray*}
%\varphi(\alpha)=\phi_{P}(\alpha)^{c},
%\end{eqnarray*}
%for all $\alpha\in {\mathbb K}$. Hence, $\varphi_{P}$ is uniquely determined by the above condition if it exists. The existence of $\varphi_P$ follows since the group $\log_{N(P)}\phi_{P}({\mathbb K})$ is infinite cyclic for any choice of $\phi_{P}\in P$ (cf. \cite{Rein}).\\
% \\
If $P\in \sum_{\mathbb K}^{\infty}$, there is a unique embedding $\theta_P:{\mathbb K}\to \widehat{{\mathbb K}}_P$ corresponding to $P$. We let $|\cdot|$ be the usual absolute value on $\widehat{{\mathbb K}}_P$ ($={\mathbb R}$ or ${\mathbb C}$), and define
\begin{eqnarray*}
\varphi_{P}=|\theta_P|^{\dim_{\mathbb R}\widehat{\mathbb K}_P}.
\end{eqnarray*}
We have the following product formula (cf. \cite{Rein}).
\bth \label{pf}
If ${\mathbb K}$ is a global field, then
\begin{eqnarray*}
\prod_{P\in\sum_{\mathbb K}}\varphi_P(\alpha)=1,
\end{eqnarray*}
for all $\alpha\in {\mathbb K}^{*}$.
\eop
\eth
Consider a global field ${\mathbb L}$, and assume that ${\mathbb K}$ is another global field such that the extension ${\mathbb L}/{\mathbb K}$ is finite and separable. To each element $Q\in\sum_{\mathbb L}$, we shall now associate an integer $e_Q$ and a real number $r_Q$, depending on this extension. \\
 \\
For $Q\in\sum_{\mathbb L}$, we denote by $P_Q$ the element in $\sum_{\mathbb K}$ that contains the restriction $\varphi_{Q}|_{\mathbb K}$. We let $B_{P_Q}$ be the integral closure of $A_{P_Q}$ in ${\mathbb L}$, and denote by $\widehat{B}_{P_Q}$ and $\widehat{A}_{P_Q}$ the corresponding completed rings.
%and denote by ${\cal D}_Q$ the different of $B_{P_Q}$ over $A_{P_Q}$ (cf. \cite{Serr}, III: \S $3$).
\newpage
\noindent
\bd \label{rami}
The {\em ramification index} of $Q$ relative to the extension ${\mathbb L}/{\mathbb K}$, is the integer
\begin{eqnarray*}
e_Q&=&\left\{ \begin{array}{lcl}
                               [\widehat{{\mathbb L}}_Q:\widehat{{\mathbb K}}_{P_Q}], &&\mbox{ if $Q\in\sum_{\mathbb L}^{\infty}$, }\\
                               (\varphi_Q({\mathbb L}^{*}):\varphi_{P_Q}({\mathbb K}^{*})), &&\mbox{ if $Q\in \sum_{\mathbb L}\setminus\sum_{\mathbb L}^{\infty}$}.
                                         \end{array} \right .
\end{eqnarray*}
For $Q\in \sum_{\mathbb L}^{\infty}$, we define
\begin{eqnarray*}
r_Q&=&\left\{ \begin{array}{lcl}
                                %\max\left(\{r\in {\mathbb Z}; (B_{P_Q}\cap M_Q)^{r}\supseteq {\cal D}_{Q}\}\right), &&\mbox{ if $Q\in\sum_{\mathbb L}\setminus\sum_{\mathbb L}^{\infty}$, }\\
                                -\log\log N(Q), &&\mbox{ if  $e_Q\neq 1$,}\\
0, &&\mbox{ otherwise}.
                                         \end{array} \right .
\end{eqnarray*}
For $Q\in\sum_{\mathbb L}\setminus\sum_{\mathbb L}^{\infty}$, we define $r_Q$ to be the exponent of $\widehat{M}_Q$ in the different of $\widehat{B}_{P_Q}$ over $\widehat{A}_{P_Q}$ (cf. {\em \cite{Serr}}, \S $3$ in Chapter III).
\eop
\ed
\bnote \label{discriminant}
When $\chara({\mathbb K})=0$, we obtain with this definition
\begin{eqnarray*}
\prod_{P\in\sum_{\mathbb K}\setminus\sum_{\mathbb K}^{\infty}}N(P)^{r_Q}=|\disc_{\mathbb K}|,
\end{eqnarray*}
where $\disc_{\mathbb K}$ denotes the discriminant of the number field ${\mathbb K}$. This follows from Proposition $6$ and Proposition $10$ in Chapter III of {\em \cite{Serr}}.
\eop
\enote
Throughout the paper, we employ the conventions of letting empty sums equal $0$, and letting empty products equal $1$.
 \\
 \\
%Claim: If ${\mathbb K}$ is a global field, then $\#\sum_{\mathbb K}^{\infty}<\infty$. (explanation)\\
\section{Divisors on global fields}\label{dogfatrr}
In this section, we define and study divisors on global fields. In particular, we consider the set of multiples of a divisor $D$, and relate its cardinality to the degree of $D$ (cf. Theorem \ref{rr}).
\newpage
\noindent
\bd \label{divdef}
Let ${\mathbb K}$ be a global field. A {\em divisor} $D$ on ${\mathbb K}$ is a formal finite sum
\begin{eqnarray*}
D=\sum_{P\in \sum_{\mathbb K}}a_P\cdot P,
\end{eqnarray*}
where $a_P\in {\mathbb R}$ if $P\in \sum_{\mathbb K}^{\infty}$, and $a_P\in {\mathbb Z}$ otherwise. The {\em degree} of $D$ is the real number
\begin{eqnarray*}
\deg D=\prod_{P\in \sum_{\mathbb K}}N(P)^{a_P}.
\end{eqnarray*}
The divisor $D$ is said to be {\em principal} if there exists an $\alpha\in{\mathbb K}^{*}$, such that
\begin{eqnarray*}
N(P)^{a_P}=\varphi_{P}(\alpha),
\end{eqnarray*}
for all $P\in\sum_{\mathbb K}$.
\eop
\ed
We may now state Theorem \ref{pf} by saying that a principal divisor has degree $1$. It follows that the identity element in the group defined by the degree homomorphism 
\begin{eqnarray*}
\deg:\{\mbox{divisors on  } {\mathbb K}\}\to {\mathbb R},
\end{eqnarray*}
is the class containing the principal divisors. The inverse of the class containing a divisor $D=\sum a_P\cdot P$, is the class containing the divisor
\begin{eqnarray*}
-D=\sum (-a_P)\cdot P.
\end{eqnarray*}
\bd \label{muldef}
Let $D=\sum a_P\cdot P$ be a divisor on a global field ${\mathbb K}$. The {\em space of multiples of $D$} is the set
\begin{eqnarray*}
H^{0}(D)=\{\alpha\in{\mathbb K}; \varphi_P(\alpha)\leq N(P)^{a_P}, \mbox{ for all $P\in\sum_{\mathbb K}$}\}.
\end{eqnarray*}
We denote by $h^{0}(D)$ the cardinality of $H^{0}(D)$.
\eop
\ed
\newpage
\noindent
\bth \label{rr}
Let $D$ be a divisor on a global field ${\mathbb K}$. \\
$(i)$ There exists a divisor $\omega_{\mathbb K}$, depending only on ${\mathbb K}$, such that
\begin{eqnarray*}
\frac{1}{C(S_{1}({\mathbb K}),S_{2}({\mathbb K}))}\leq \frac{h^{0}(D)}{h^{0}(\omega_{\mathbb K}-D)}\cdot \frac{\sqrt{\deg \omega_{\mathbb K}}}{\deg D}\leq C(S_{1}({\mathbb K}),S_{2}({\mathbb K})),
\end{eqnarray*}
with
\begin{eqnarray*}
C(S_{1}({\mathbb K}),S_{2}({\mathbb K}))=\frac{6^{S_{1}({\mathbb K})+S_{2}({\mathbb K})}\cdot (S_{1}({\mathbb K})+S_{2}({\mathbb K}))!}{2^{S_{1}({\mathbb K})}\cdot (\pi/2)^{S_{2}({\mathbb K})}}.
\end{eqnarray*}
$(ii)$ There exists a function $i:$ \{divisors on ${\mathbb K}\}\to {\mathbb R}$, such that $i(\cdot)\to 1$ when $\deg \cdot\to \infty$, and
\begin{eqnarray*}
\frac{h^{0}(D)}{i(D)}\cdot\frac{\sqrt{\deg \omega_{\mathbb K}}}{\deg D}=2^{S_{1}({\mathbb K})}\cdot (2\pi)^{S_{2}({\mathbb K})/2}.
\end{eqnarray*}
\eth
\prf 
Assume first that $\chara({\mathbb K})=0$, and denote by ${\cal O}_{\mathbb K}$ the integral closure of ${\mathbb Z}$ in ${\mathbb K}$. Consider the ${\cal O}_{\mathbb K}$-module $\homo_{\mathbb Z}({\cal O}_{\mathbb K}, {\mathbb Z})$, metrized by defining \\$|\trace|_P=\log N(P)$, for $P\in \sum_{\mathbb K}^{\infty}$ (cf. \cite{GiSo}, \S$2.4$). \\
If a divisor $\omega_{\mathbb K}$ on ${\mathbb K}$ is chosen such that 
\begin{eqnarray*}
\deg\omega_{\mathbb K}=\frac{|\disc_{\mathbb K}|}{2^{S_{2}({\mathbb K})}},
\end{eqnarray*}
the corresponding metrized ${\cal O}_{\mathbb K}$-module will be isometrically isomorphic to \\$\homo_{\mathbb Z}({\cal O}_{\mathbb K}, {\mathbb Z})$, metrized as above (cf. \cite{Neu}, Theorem $4.5$ in Chapter III). \\Hence one obtains $(i)$ from Theorem $2$ in \cite{GiSo}. However, note Remark \ref{correction} on the value of $C(S_{1}({\mathbb K}),S_{2}({\mathbb K}))$.\\
 \\
For a divisor $D=\sum a_P\cdot P$, denote by $\chi(D)$ the Euler-Minkowski characteristic (cf. \cite{Neu}, Definition $3.1$ in \S $3$ of Chapter III) of the fractional ideal
\begin{eqnarray*}
\prod_{P\in\sum_{\mathbb K}\setminus\sum_{\mathbb K}^{\infty}}({\cal O}_{\mathbb K}\cap M_{P})^{-a_P}.
\end{eqnarray*}
Setting
\begin{eqnarray*}
i(D)=\frac{h^{0}(D)\cdot e^{-\chi(D)}}{2^{S_{1}({\mathbb K})}(\pi)^{S_{2}({\mathbb K})/2}},
\end{eqnarray*}
one obtains $(ii)$ as a slight reformulation of Theorem $3.9$ (Chapter III, \S $3$) in \cite{Neu}.\\
 \\
Now assume that $\chara({\mathbb K})>0$. In this case $S_{1}({\mathbb K})=S_{2}({\mathbb K})=0$, and
\begin{eqnarray*}
\log_q h^{0}(D)=\dim_{{\mathbb F}_q}H^{0}(D). 
\end{eqnarray*}
\newpage
\noindent
Let $g_{\mathbb K}$ be the genus of the complete non-singular curve determined by ${\mathbb K}$, and choose $\omega_{\mathbb K}$ from the class of divisors of degree $q^{2g_{\mathbb K}-2}$. Then $(i)$ is a multiplicative formulation of the Riemann-Roch theorem (cf. \cite{Rose}, Theorem $5.4$ in Chapter $5$). Setting $i(D)=h^{0}(\omega_{\mathbb K}-D)$, $(ii)$ follows from $(i)$, since $h^{0}(\omega_{\mathbb K}-D)=1$ whenever $\deg D>\deg\omega_{\mathbb K}$ (cf. \cite{Rose}, Corollary $4$ in Chapter $5$).
\eop
\bnote \label{correction}
We make a minor correction to the proof of Theorem $2$ in {\em \cite{GiSo}}. Numbers in bold-face refer to lines or pages in {\em \cite{GiSo}}. The quantity $C(r_1 , r_2 , N)$ is defined on line {\em (}{\bf 26}{\em )} {\em (}pg. {\bf 355}{\em )} as
\begin{eqnarray*}
-\log\mu(K^{*})+N (r_1+2r_2) \log (6),
\end{eqnarray*}
where $\mu(K^{*})$ is the euclidean volume of the set of  $({\bf y}_i,{\bf z}_j)\in ({\mathbb R}^{N})^{r_1}\times ({\mathbb C}^{N})^{r_2}$ such that
\begin{eqnarray*}
\sum_{i=1}^{r_1}|{\bf y}_i| + 2\sum_{j=1}^{r_2}|{\bf z}_j|\leq 1.
\end{eqnarray*}
However, the value of $C(r_1,r_2,N)$ stated in Theorem $2$ in {\em \cite{GiSo}} is incorrect, due to a missing minus sign in the computation of $\mu(K^{*})$ on line {\em (}{\bf 24}{\em )} {\em (}pg. {\bf 355}{\em )}. The correct value is
\begin{eqnarray*}
C(r_1,r_2,N)&=&\log\left(\frac{6^{N(r_1+2r_2)}}{\mu(K^{*})}\right)\\
&=&\log\left(\frac{(N(r_1+2r_2))!\cdot 2^{2Nr_2}\cdot 6^{N(r_1+2r_2)}}{(V(B_N) N!)^{r_1}\cdot (V(B_{2N}) (2N)!)^{r_2}}\right).
\end{eqnarray*}
The value of $C(S_{1}({\mathbb K}),S_{2}({\mathbb K}))$ in Theorem \ref{rr} is simply $e^{C(r_1,r_2,1)}$. 
\eop
\enote
\section{A canonical divisor}\label{atrht}
In this section, we describe a divisor $\omega_{\mathbb K}'$ that is determined by the global field ${\mathbb K}$, and show that Theorem \ref{rr} holds with $\omega_{\mathbb K}=\omega_{\mathbb K}'$. We also show that for a finite separable extension ${\mathbb L}/{\mathbb K}$ of global fields, the corresponding divisors $\omega_{{\mathbb L}}'$ and $\omega_{{\mathbb K}}'$ satisfy a Riemann-Hurwitz type formula (cf. Theorem \ref{rh}).\\
 \\
 We begin by considering a divisor that is determined by a finite separable extension of global fields. Recall the real numbers $r_Q$ from Definition \ref{rami}.\newpage
\noindent
\bd
Let ${\mathbb L}/{\mathbb K}$ be a finite separable extension of global fields. The {\em ramification divisor} relative to the extension ${\mathbb L}/{\mathbb K}$ is the divisor
\begin{eqnarray*}
R_{{\mathbb L}/{\mathbb K}}=\sum_{Q\in {\mathbb L}}r_Q\cdot Q.
\end{eqnarray*}
\eop
\ed
If $\chara({\mathbb K})>0$, we denote by $P_0$ a fixed element of $\sum_{{\mathbb K}_0}\setminus\sum_{{\mathbb K}_0}^{\infty}$ such that $N(P_0)=q$. If $\chara({\mathbb K})=0$, we let $P_0$ be an arbitrary fixed element of $\sum_{{\mathbb K}_0}\setminus\sum_{{\mathbb K}_0}^{\infty}$. In both cases, we denote by $S_0$ the set of $P\in\sum_{\mathbb K}$ such that the restriction $\varphi_P|_{{\mathbb K}_0}$ is contained in $P_0$.\\
 \\
If  $\chara({\mathbb K})=0$, we choose in addition an element $P_{\infty}\in\sum_{\mathbb K}^{\infty}$, and set
\begin{eqnarray*}
a_{\infty} =\sum_{P\in S_0}2e_P\log_{N(P_{\infty})}N(P),
\end{eqnarray*}
where $e_P$ denotes the ramification index of $P$ relative to the extension ${\mathbb K}/{\mathbb K}_0$ (cf. Definition \ref{rami}).\\
 \\
Recall that the extension ${\mathbb K}/{\mathbb K}_0$ is separable by definition (cf. Remark \ref{knoll}), and consider the divisor
\begin{eqnarray*}
\omega_{\mathbb K}'=R_{{\mathbb K}/{\mathbb K}_0}-\sum_{P\in S_0}2e_P\cdot P + a_{\infty}\cdot P_{\infty}.
\end{eqnarray*}
\bpr Theorem \ref{rr} holds with $\omega_{\mathbb K}=\omega_{\mathbb K}'$.\epr
\prf It suffices to verify that
\begin{eqnarray*}
\deg \omega_{\mathbb K}'&=&
\left\{ \begin{array}{lcl}
                                |\disc_{\mathbb K}|\cdot 2^{-S_{2}({\mathbb K})}, &&\mbox{ if $\chara(\mathbb K)=0$, }\\
                               q^{2g_{\mathbb K}-2}, &&\mbox{ if $\chara({\mathbb K})>0$}.
  \end{array} \right .
\end{eqnarray*}
%(Here $\disc_{\mathbb K}$ denotes the discriminant of ${\mathbb K}$ over ${\mathbb Q}$ when $\chara({\mathbb K})=0$, and $g_{\mathbb K}$ denotes the genus of the complete non-singular curve determined by ${\mathbb K}$ when $\chara({\mathbb K})>0$.)\\
%Recall that : for $Q\in\sum_{\mathbb K}\setminus\sum_{\mathbb K}^{\infty}$, denote by $P_Q$ the element of $\sum_{{\mathbb K}_0}$ containing the restriction $\varphi_{Q}|_{{\mathbb K}_0}$. Denote by $B_{P_Q}$ the integral closure of $A_{P_Q}$ in ${\mathbb K}$, and by ${\cal D}_{B_{P_Q}/A_{P_Q}}$ the different of $B_{P_Q}$ over $A_{P_Q}$. 
If $\chara({\mathbb K})=0$, one has
\begin{eqnarray*}
\deg{\omega_{\mathbb K}'}=e^{-S_{2}({\mathbb K})\log 2}\prod_{Q\in\sum_{\mathbb K}\setminus\sum_{\mathbb K}^{\infty}}N(Q)^{r_Q}=\frac{|\disc_{\mathbb K}|}{2^{S_{2}({\mathbb K})}},
\end{eqnarray*}
where Remark \ref{discriminant} is used to obtain the last equality.\\
 \\
If $\chara({\mathbb K})>0$, one has $g_{{\mathbb K}_0}=0$. Hence
\begin{eqnarray*}
\deg{R_{{\mathbb K}/{\mathbb K}_0}}=q^{2g_{\mathbb K}-2+2[{\mathbb K}:{\mathbb K}_0]},
\end{eqnarray*}
by the Riemann-Hurwitz formula for function fields (cf. \cite{Rose}, Theorem $7.16$ in Chapter $7$). \\
\newpage
\noindent
Since the extension ${\mathbb K}/{\mathbb K}_0$ is separable, one has
\begin{eqnarray*}
\prod_{P|P_0} N(P)^{e_P}=q^{[{\mathbb K}:{\mathbb K}_0]},
\end{eqnarray*}
for any choice of $P_0\in {\mathbb K}_0$ such that $N(P_0)=q$ (cf. \cite{Serr}, Proposition $10$ in \S$4$ of Chapter I). This completes the proof.
\eop
%... , we obtain the following extension of the Riemann-Hurwitz formula (cf. ..).
\bth \label{rh}
If ${\mathbb L}/{\mathbb K}$ is a finite separable extension of global fields, and if ${\mathbb L}_0={\mathbb K}_0$, then
\begin{eqnarray*}
\deg{\omega_{\mathbb L}'}=\deg{\omega_{\mathbb K}'}^{[{\mathbb L}:{\mathbb K}]}\cdot\deg{R_{{\mathbb L}/{\mathbb K}}}.%\cdot \sqrt{2^{S_{\infty}(2)([{\mathbb L}:{\mathbb K}]-1)}}.
\end{eqnarray*}
\eth
\prf
Assume first that $\chara({\mathbb K})=0$. Denote the restrictions of $R_{{\mathbb L}/{\mathbb K}}$, $\omega_{\mathbb L}'$ and $\omega_{\mathbb K}'$ to the archimedean classes by $R_{{\mathbb L}/{\mathbb K}}^{\infty}$, $\omega_{\mathbb L}^{\infty}$ and $\omega_{\mathbb K}^{\infty}$, respectively. Note that by the construction of $R_{{\mathbb L}/{\mathbb K}}$, $\omega_{\mathbb L}'$ and $\omega_{\mathbb K}'$:
\begin{eqnarray*}
(i) &&\mbox{$\log_{1/4}\deg R_{{\mathbb L}/{\mathbb K}}^{\infty}$ is the number of elements $Q\in\sum_{\mathbb L}^{\infty}$ with}\\ 
&&\mbox{$\log N(Q)=2$ extending elements $P\in\sum_{\mathbb K}^{\infty}$ with $\log N(P)=1$,}\\
(ii) &&\mbox{${[{\mathbb L}:{\mathbb K}]}\cdot\log_{1/4}\deg \omega_{\mathbb K}^{\infty}$ is the number of elements $Q\in\sum_{\mathbb L}^{\infty}$ with}\\
&&\mbox{$\log N(Q)=2$ extending elements $P\in\sum_{\mathbb K}^{\infty}$ with $\log N(P)=2$,}\\
(iii) &&\mbox{$\log_{1/4}\deg \omega_{\mathbb L}^{\infty}$ is the total number of elements $Q\in\sum_{\mathbb L}^{\infty}$ with}\\&&\mbox{$\log N(Q)=2$.} 
\end{eqnarray*}
From these remarks, we obtain the equality
\begin{eqnarray*}
\deg{\omega_{\mathbb L}^{\infty}}=\deg{\omega_{\mathbb K}^{\infty}}^{[{\mathbb L}:{\mathbb K}]}\cdot\deg{R_{{\mathbb L}/{\mathbb K}}^{\infty}}.
\end{eqnarray*}
The corresponding equality for $R_{{\mathbb L}/{\mathbb K}}-R_{{\mathbb L}/{\mathbb K}}^{\infty}$, $\omega_{\mathbb L}'-\omega_{\mathbb L}^{\infty}$ and $\omega_{\mathbb K}'-\omega_{\mathbb K}^{\infty}$ follows from the transitivity of the different in a tower of finite separable extensions of fields (cf. \cite{Serr}, Proposition $8$ in \S $4$ of Chapter III).\\
 \\
When $\chara({\mathbb K})>0$, the statement in the theorem is a multiplicative formulation of the Riemann-Hurwitz formula for function fields (cf. \cite{Rose}, Theorem $7.16$ in Chapter $7$).
\eop
\newpage
\noindent
\bibliography{bibl}
\end{document}